\documentclass[12
pt]{article}

\usepackage[latin1]{inputenc}

 \usepackage{graphicx} 
\usepackage{amsfonts,amsmath,amssymb,amsthm,graphicx,epsfig,float}
\usepackage[T1]{fontenc}

\usepackage[francais]{babel}

{
 \newtheorem{theoreme}{Th\'eor\`eme}
  \newtheorem*{theoreme*}{Th\'eor\`eme}
  \newtheorem{lemme}[theoreme]{Lemme}
  \newtheorem{corollaire}[theoreme]{Corollaire}
  
 \newtheorem{conjecture}[theoreme]{Conjecture}
  \newtheorem{definition}{D\'efinition}
  \newtheorem{proposition}[theoreme]{Proposition}
\theoremstyle{remark}
  \newtheorem*{remarque*}{Remarque}
}

\newcommand{\Nn}{\mathbb{N}}

\renewcommand {\epsilon}{\varepsilon}

\renewcommand {\leq}{\leqslant}
\renewcommand {\geq}{\geqslant}

\title{{\bf A propos de la conjecture de Nash }}
\author{Camille Plénat}
\date{\today}

\begin{document}
\maketitle
\begin{abstract}

Cet article a pour objet le probl\`eme des arcs de Nash, selon lequel il y aurait autant de familles d'arcs sur un germe de surface singulier $U$ que de composantes essentielles d'une d\'esingularisation de cette singularit\'e. Soit $\mathcal{H}=\bigcup \overline{N_\alpha}$ une d\'ecomposition particuli\`ere de l'espace des arcs, d\'ecrite plus avant. On propose deux nouvelles conditions suffisantes permettant d'affirmer que $\overline{N_\alpha}\not \subset \overline{N_\beta}$, $\alpha \not = \beta$; plus pr\'ecisemment,pour la premi\`ere, on montre que, s'il existe une fonction f dans l'anneau local de $U$ telle que $ord_{E_\alpha}(f)<ord_{E_\beta}(f)$, o\`u $E_\alpha,\ E_\beta$ sont des diviseurs issus de la d\'esingularisation minimale, alors $\overline{N_\alpha}\not \subset \overline{N_\beta}$.La deuxi\`eme condition, utilis\'ee dans le cadre d'une singularit\'e rationnelle de surface, est la suivante: soit $(S,s)$ et $(S',s')$ deux singularit\'es rationnelles de surface telles qu'il existe un morphisme dominant et birationnel $\pi$ de $(S,s)$ sur $(S',s')$;alors, soient $E_\alpha,\ E_\beta$ deux diviseurs issus de la d\'esingularisation minimale de $(S,s)$ tels que leurs images par $\pi$, $E'_\alpha$ et $E'_\beta$, sont des diviseurs exceptionnels pour $(S',s')$; si $\overline{N'_\alpha}(S',s')\not \subset \overline{N'_\beta}(S',s')$ alors $\overline{N_\alpha}(S,s)\not \subset \overline{N_\beta}(S,s)$.  Ces conditions, certes tr\`es simples, nous permettent de prouver assez directement la conjecture pour les singularit\'es rationnelles minimales de surface (en utilisant conjointement la d\'ecomposition des singularit\'es rationnelles minimales de surface en singularit\'es \`a quotient cyclique de type $A_n$). Une preuve de la conjecture pour ces singularit\'es a d\'ej\`a \'et\'e donn\'ee par Ana Reguera dans [{\bf REG}].
\end{abstract}
\def\figurename{{Fig.}}%
\def\contentsname{Sommaire}%

%
{\bf Introduction}
\hspace{2cm}

L'objet de ce texte est d'\'enoncer deux nouvelles conditions suffisantes pour pouvoir affirmer que les adh\'erences de deux familles d'arcs sur une singularit\'e associ\'ees \`a deux diviseurs exceptionnels d'une d\'esingularisation ne sont pas incluses l'une dans l'autre. Pour les singularit\'es rationnelles de surface, la premi\`ere condition que je vais \'enoncer est \'equivalente \`a celle donn\'ee par Ana Reguera dans [{\bf REG}]. Comme applications directes, nous allons tout d'abord d\'efinir une relation binaire sur les diviseurs exceptionnels (qui sera un ordre partiel pour les singularit\'es rationnelles) et donc identifier les familles ``\`a probl\`emes'', puis nous allons donner une preuve tr\`es simple qu'une r\'eponse affirmative peut \^etre donn\'ee au probl\`eme de Nash pour les singularit\'es minimales de surface, r\'esultat d\'ej\`a acquis par Ana Reguera (cf{\bf REG}): dans un premier temps nous d\'emontrons la conjecture pour les singularit\'es quotient cyclique de type $A_n$ en appliquant la premi\`ere condition, puis nous recourons \`a la d\'ecomposition des singularit\'es minimales en singularit\'es quotient cyclique de type $A_n$  (cf. [{\bf SPI 1}])ainsi qu'\`a la deuxi\`eme condition, afin de prouver la conjecture  pour ces singularit\'es.\\\\
{\bf Remarque:}la d\'ecomposition des singularit\'es minimales de surface en singularit\'es quotient cyclique de type $A_n$ est un cas particulier de la d\'ecomposition des singularit\'es sandwichs en singularit\'es primitives. Cette nouvelle d\'emonstration de la conjecture pour les singularit\'es minimales nous donne un espoir de prouver un jour la conjecture pour les singularit\'es sandwichs, sous r\'eserve de la prouver pour les primitives.\\\\ 
Je tiens \`a remercier ici mon directeur de th\`ese, Mark Spivakovsky, qui m'a constamment guid\'ee dans l'\'elaboration de cet article (et mes camarades du laboratoire Emile Picard dont l'aide a \'et\'e pr\'ecieuse pour sa mise en forme).

\section{{\bf Rappels et definitions}}

Dans cette section on rappellera les d\'efinitions ayant trait \`a la conjecture de Nash ainsi que la conjecture elle--m\^eme.\\

\subsection{{\bf La conjecture de Nash}}

Cette conjecture concerne la famille des courbes formelles passant par une singularit\'e et les composantes exceptionnelles d'une d\'esingularisation de cette singularit\'e. La conjecture a \'et\'e \'enonc\'ee en toute g\'en\'eralit\'e, bien que nous ne l'\'enoncerons ici que pour les singularit\'es de surfaces.\\
Avant d'\'enoncer la conjecture, rappelons quelques d\'efinitions:\\

\begin{definition}{\bf Arc}\\
Un arc sur un germe $U$ est une courbe formelle param\'etris\'ee, i.e. un k--morphisme $\phi:(spec\ k[[t]],0)\longrightarrow(S,U)$, ou, de mani\`ere \'equivalente, un k--morphisme ${\cal O}_{U}\longrightarrow k[[t]]$ .\\
\end{definition}

Soit $(S,s)$ une singularit\'e isol\'ee de surface.\\

\begin{definition}
On d\'efinit l'image d'un arc $\phi$ comme \'etant $Z(\phi)=sC$ o\`u $C$ est le diviseur de Weil d\'efini par $\phi(spec\ k[[t]])$, et si $\overline{C}$ est la normalisation de $C$, alors $s$ est l'indice de ramification du morphisme $(spec\ k[[t]],0)\longrightarrow$$(\overline{C},P)$.\\
\end{definition}

\begin{definition}
Soit $\mathcal{H}$ l'ensemble des arcs du germe singulier $U$. C'est une vari\'et\'e alg\'ebrique de dimension infinie.\\
Quels sont les ouverts sur $\mathcal{H}$? Ils sont d\'efinis par un nombre fini d'in\'equations sur un nombre fini de coefficients des s\'eries enti\`eres formelles d\'efinissant les arcs.\\
\end{definition}
(cf.infra section 2.1 pour un calcul explicite de  $\mathcal{H}$ pour les singularit\'es $A_n$).\\\\
Soit de nouveau $(S,s)$ une singularit\'e.
Soient $(X_0,{\{E_{\alpha}\}}_{(\alpha)})\longrightarrow(S,s)$ une d\'esingularisation de la singularit\'e $(S,s)$, et $\{E_{\alpha}\}_{(\alpha)}$ l'ensemble des composantes exceptionnelles associ\'ees \`a la d\'esingularisation.\\
La conjecture de Nash consiste \`a comparer les composantes irr\'eductibles de  $\mathcal{H}$ et les composantes $\{E_{\alpha}\}_{(\alpha)}$; plus exactement, dans le cas des surfaces:
\begin{conjecture}{\bf (Nash)}\\
Il y a autant de courbes  exceptionnelles dans la {\bf d\'esingularisation minimale} d'une singularit\'e de surface $(S,s)$ que de composantes irr\'eductibles de $\mathcal{H}$.
\end{conjecture}

Ces composantes irr\'eductibles sont aussi appel\'ees familles d'arcs de Nash.
En associant \`a chaque famille d'arcs un diviseur exceptionnel, Nash, dans [{\bf NASH}], a montr\'e que le nombre de familles est major\'ee par le nombre de composantes exceptionnelles.\\
 Consid\'erons la d\'ecomposition de  $\mathcal{H}$ suivante:\\
soit $N_{\alpha}$ la famille d'arcs tels que leur transform\'ee stricte est transverse \`a $E_{\alpha}$ et n'intersecte pas les autres composantes exceptionnelles ; ces familles sont des ensembles irr\'eductibles de $\mathcal{H}$ et  $\mathcal{H}=\bigcup \overline{N_{\alpha}}$ ;remarquons qu'il y a d\'ej\`a dans cette d\'ecomposition autant de familles que de diviseurs, il nous suffit donc de prouver que pour tout $\alpha ,\gamma $, on a $\overline{N_{\alpha}}\not\subset \overline{N_{\gamma}}$.\\
Pour de plus amples explications, voir l'article d'Ana Reguera [{\bf REG}].

\subsection{\bf Quelques rappels sur les singularit\'es}

\hspace{2cm}Pour plus de coh\'erence, on rappellera dans ce paragraphe les d\'efinitions des singularit\'es primitives, minimales, quotient cyclique de type $A_n$, sandwichs, et la d\'ecomposition de ces derni\`eres en joint birationnel ; pour une th\'eorie plus compl\`ete ainsi que pour toutes les d\'efinitions, on renvoie le lecteur aux articles de Mark Spivakovsky (cf.[{\bf SPI 1},premi\`ere partie] et [{\bf SPI 2}]). De plus on ne donne les d\'efinitions qu'en ``fonction'' des graphes car les d\'emonstrations sont esentiellement bas\'ees sur la nature de ces graphes. \\

{\bf Rappel:} quand on d\'esingularise une singularit\'e rationnelle  normale de surface, des courbes exceptionnelles apparaissent et forment un graphe connexe sans cycle. Mais on n'\'etudie pas ces graphes l\`a, on pr\'ef\`ere regarder les graphes duaux, c'est--\`a--dire ceux qui associent \`a chaque point une courbe exceptionnelle, chaque ar\^ete repr\'esentant alors l'intersection entre deux courbes.\\ Les poids que l'on associera \`a une courbe sont moins le nombre d'autointersection de cette courbe avec elle--m\^eme.(voir fig.1)

\begin{figure}[H]
$$
\setlength{\unitlength}{1cm}
\begin{picture}(8,4)
\includegraphics{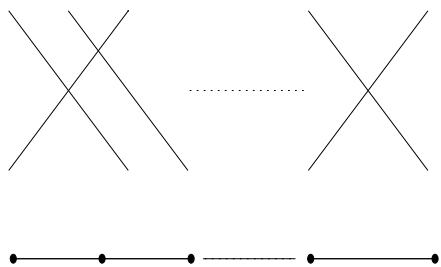}
\end{picture}
$$
\caption{exemple: en haut, graphe avec les composantes; en dessous, graphe dual associ\'e (singularit\'e de type $A_n$)}
\end{figure}

{\bf Remarque}: comme on ne va consid\'erer que les d\'esingularisations minimales, les poids seront sup\'erieurs ou \'egaux \`a 2 (sauf pour les graphes non singuliers, i.e. issus de la d\'esingularisation d'un point non singulier).\\ 


\begin{definition}{\bf Graphe sandwich}\\
Un graphe est dit sandwich s'il existe un graphe non singulier le contenant comme sous--graphe pond\'er\'e ;
une singularit\'e sera dite sandwich si son graphe dual de r\'esolution est sandwich.\\
\end{definition}

\begin{definition}{\bf Graphe primitif}\\
Un graphe $\Gamma$ connexe est dit primitif s'il existe un graphe non singulier ${\Gamma}^{*}$ tel que $\Gamma \subset {\Gamma}^{*}$ et $\#({\Gamma}^{*} \backslash \Gamma )$$=1$.\\
Une singularit\'e sera dite primitive si son graphe dual de r\'esolution est primitif.\\
\end{definition}

\begin{definition}{\bf Graphe minimal}\\
Un graphe pond\'er\'e est dit minimal  si pour tout x sommet du graphe  $\omega (x) \geq \gamma (x)$, o\`u $\omega (x)$ est le poids en $x$ et $ \gamma (x)$ est le nombre d'ar\^etes attach\'ees \`a x.\\
 Une singularit\'e est dite minimale si son graphe dual de r\'esolution est minimal, en d'autres termes la singularit\'e est minimale si dans sa d\'esingularisation minimale toutes les courbes exceptionnelles sont rationnelles, s'intersectent transversalement et le graphe dual est simplement connexe et minimal.\\

{\bf Remarque:} cela \'equivaut \`a dire que la singularit\'e est rationnelle et que son cycle fondamental est r\'eduit.\\
\end{definition}

\begin{definition}{\bf Singularit\'e quotient cyclique de type $A_n$}\\
Par $A_n$ on d\'esigne les singularit\'es quotient cyclique dont les poids \`a chaque sommet du graphe dual valent 2 . Ce sont des singularit\'es de surface d'\'equations ($z^{n+1}=xy$) ; leur graphe dual de r\'esolution est:

\begin{figure}[H]
$$
\setlength{\unitlength}{1cm}
\begin{picture}(12,1)
\includegraphics{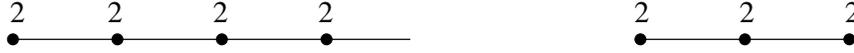}
\end{picture}
$$
\caption{graphe dual d'une singularit\'e de type $A_n$}
\end{figure}

On appellera ce genre de graphe un graphe ``en bambou''.
\end{definition}

\begin{definition}{\bf Joint birationnel}\\
Soit $Z_1 ,Z_2$ deux vari\'etes alg\'ebriques irr\'eductibles et r\'eduites avec le m\^eme corps de fonction K.
Soit $\phi:Z_1 \leftrightarrow Z_2 $ la correspondance birationnelle entre ces deux vari\'et\'es. Alors il existe deux ouverts $U_1$ et $U_2$ tels que $\phi$ induit un isomorphisme entre $U_1$ et $U_2$. Soit U l' image de $U_1$ dans $U_1\times U_2$ par le morphisme diagonal o\`u on identifie $U_1$ avec $U_2$ par $\phi$. Alors le joint birationnel Z de $Z_1$ avec $Z_2$ est d\'efini comme \'etant la fermeture de U dans $Z_1\times Z_2$.\\
\end{definition}

Les singularit\'es minimales sont sandwichs, les quotients cycliques de type $A_n$ sont des primitives, lesquelles sont elles--m\^emes sandwichs. On peut regarder les singularit\'es sandwichs (resp. minimales) comme le joint birationnel de singularit\'es primitives (resp. \`a quotient cyclique de type $A_n$).\\\\
{\bf Remarque:} Dans les faits on regardera les d\'ecompositions en joint birationnel des singularit\'es sandwichs en singularit\'es primitives et des minimales en ``$A_n$'' au niveau des graphes. Nous verrons dans la preuve de la conjecture pour les singularit\'es minimales comment dans la pratique on peut d\'ecomposer les graphes des minimales en graphe de singularit\'es $A_n$.
 
\begin{lemme}
Les singularit\'es $A_n$ sont des  singularit\'es sandwichs.
\end{lemme}

\begin{proof}
Le graphe de la singularit\'e $A_n$ est contenu dans un graphe non singulier, comme le montre la figure ci--dessous:
(le signe $<$ veut dire ici ``est contenu dans'')
\begin{figure}[H]
$$
\setlength{\unitlength}{1cm}
\begin{picture}(12,1)
\includegraphics{exemple3.pstex}
\end{picture}
$$
\end{figure}

Donc par d\'efinition d'une singularit\'e sandwich, une singularit\'e de type $A_n$ est une singularit\'e sandwich.

\end{proof}

\section {{\bf Premi\`ere condition suffisante et r\'esolution de la conjecture pour les singularit\'es de type $A_n$}} 
Dans un premier temps, nous \'enoncerons une condition suffisante annonc\'ee pour que  $\overline{N_i}(S)\not\subset\overline{N_j}(S)$ pour une singularit\'e quelconque, puis on appliquera cette condition aux singularit\'es de type $A_n$.   

\subsection {{\bf Condition suffisante}}

\begin{definition}
Soit $S$ une vari\'et\'e singuli\`ere, $U\subset S$ ouvert affine. Soient $p:X \rightarrow S$ une d\'esingularisation de $S$ et $E_i$ un diviseur issu de la d\'esingularisation de $S$ tel que $p^{-1}(U)\cap E_i$ est non vide.\\
Soit  $ f\in {\mathcal{O}}_{U}$. On d\'efinit l'ordre de $f$ sur le diviseur $E_i$, not\'e $ord_{E_i}(f)$, comme suit: si $f^{*}=p^{*}(f) ={\mu}_i E_i+(f^{'})$ (avec $ (f^{'})$ n'ayant pas $E_i$ comme composante) on note ${\mu}_i=ord_{E_i}(f)$.\\
\end{definition}

\begin{lemme}
Soit $S$ une vari\'et\'e singuli\`ere, $U\subset S$ ouvert affine. Soient $p:X \rightarrow S$ une d\'esingularisation de $S$, $E_i$, $E_j$ diviseurs issus de la d\'esingularisation de $S$ tels que $p^{-1}(U)\cap E_k$ est non vide.\\
S'il existe $ f\in {\mathcal{O}}_{U}$ telle que $ord_{E_i}(f)< ord_{E_j}(f)$, alors  $\overline{N_i}(S)\not\subset\overline{N_j}(S)$.
\end{lemme}

\begin{proof}
Soit  $ f\in {\mathcal{O}}_{U}$ telle que $ord_{E_i}(f)< ord_{E_j}(f)$.
Alors, si $\phi$ est un arc, $ord_t (\phi(f))=Z(\phi).(f)= Z(\phi)^{'}.(f^{*}) $ (cf. [{\bf REG}]) (on d\'esigne par $Z(\phi)^{'} $ la transform\'ee stricte de $Z(\phi) $  par $p$).\\
Donc, si on note $f^{*}=p^{*}(f) ={\mu}_i E_i+(f^{'})$ comme pr\'ec\'edemment, $ord_t (\phi(f))={\mu}_i E_i.Z(\phi)^{'}+(f^{'}).Z(\phi)^{'}$. Prenons $\phi={\phi}_i\in N_i$ ne rencontrant pas $(f')$, alors $(f^{'}).Z(\phi)^{'}=0$ et ${\mu}_i E_i.Z(\phi)^{'}={\mu}_i$, donc finalement ${\bf ord_t ({\phi}(f))={\mu}_i=ord_{E_i}(f)}$.\\
Par cons\'equent si  $ord_{E_i}(f)< ord_{E_j}(f)$, pour tout ${\phi}_i\in N_i$ et ${\phi}_j\in N_j$ ne rencontrant pas $(f')$ 

\begin{eqnarray}
\ ord_{t}({\phi}_i(f))< ord_{t}({\phi}_j(f))
\end{eqnarray}

 Les coefficients de ${\phi}_k (f)$ sont des polyn\^omes en les coefficients de ${\phi}_k$. Donc l'in\'egalit\'e $(1)$ donne une in\'equation sur les coefficients de  ${\phi}_i$ que les coefficients de ${\phi}_j$ ne v\'erifient pas (le ${\mu}_i\`eme$ coefficient est non nul). Alors $\overline{N_i}\not \subset\overline{N_j}$.\\
En effet si $\overline{N_i} \subset\overline{N_j}$, prenons  ${\phi}_i\in N_i$, alors, pour tout ouvert voisinage $V_k$ de ${\phi}_i$ $V_k\cap N_j\not = \emptyset$. Soit  $V_k$ d\'efini par l' in\'equation ci--dessus, c'est un voisinage de ${\phi}_i$ mais il n'intersecte pas $N_j$.$=><=$.
\end{proof}

\subsection {{\bf Premi\`ere application: instauration d'une relation binaire sur les diviseurs exceptionnels}}
Le crit\`ere pr\'ec\'edent nous permet d'instaurer une relation binaire transitive sur l'ensemble des diviseurs exceptionnels de la d\'esingularisation minimale d'une singularit\'e de surface quelconque (qui sera un ordre partiel pour les singularit\'es rationnelles). Cela se fait de la mani\`ere suivante:\\
soient $E_i$ et $E_j$ deux diviseurs exceptionnels. Alors:\\
\begin{itemize}
   \item s'il existe f et g dans l'anneau local de la singularit\'e telles que  $ord_{E_i}(f)< ord_{E_j}(f)$ et  $ord_{E_j}(g)< ord_{E_i}(g)$, on dira que  $E_i$ et $E_j$ sont {\bf incomparables};
   \item sinon pour tout f de l'anneau local on a  $ord_{E_i}(f)\leq ord_{E_j}(f)$. Si les ordres sont \'egaux pour tout f, le crit\`ere suffisant ne donne aucune information. S'il existe un f telle que l'in\'egalit\'e est stricte, alors on dira que $E_i<E_j$\\
\end{itemize}
{\bf Remarque}:Dans le cas des singularit\'es rationnelles, il est impossible d'avoir \'egalit\'e des ordres $ord_{E_i}(f)$ et $ ord_{E_j}(f)$ pour tout f. En effet, la matrice d'intersection est d\'efinie n\'egative, donc non singuli\`ere. Donc il existe un cycle exceptionnel ${\sum}_im_iE_i=B$ tel que $B.E_k\leq0$ pour tout k et $m_i \not=m_j$. Par le th\'eor\`eme de Artin, il existe f telle que $ord_{E_i}(f)=m_i$ et $ ord_{E_j}(f)=m_j$.\\\\
Cela nous permet d'identifier facilement (on est ramen\'e \`a de l'alg\`ebre lin\'eaire) les familles de courbes pour lesquelles la conjecture est v\'erif\'ee (et donc aussi celles pour lesquelles la conjecture n'est pas encore v\'erif\'ee); en effet:\\
- si $E_i$ et $E_j$ sont incomparables alors $\overline{N_i} \not\subset\overline{N_j}$ et $\overline{N_j} \not\subset\overline{N_i}$;\\
- si $E_i<E_j$ alors $\overline{N_i} \not\subset\overline{N_j}$; mais on ne peut conclure pour la r\'eciproque.\\ 
- si $E_i=E_j$, on ne peut conclure dans aucun sens.\\

On peut sch\'ematiser cet ordre partiel sur un graphe dont nous allons donner un premier exemple, la singularit\'e $E_6$ donn\'ee par l'\'equation $z^2+y^3+x^4=0$:\\

\begin{figure}[H]
$$
\setlength{\unitlength}{1cm}
\begin{picture}(14,6)
\includegraphics{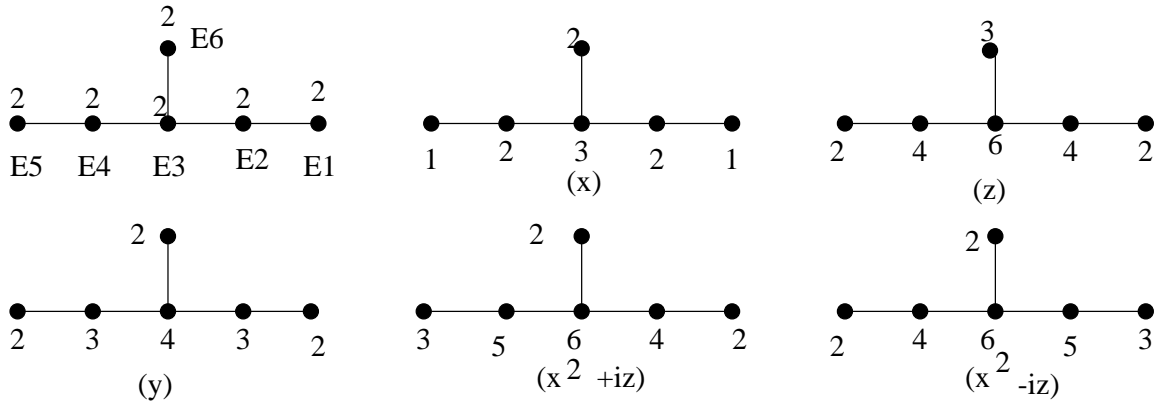}
\end{picture}
$$
\caption{Graphe dual de $E_6$ et repr\'esentations des diviseurs $(x),(y),(Z),(x^2+iz),(x^2-iz)$}
\end{figure}

En appliquant la condition suffisante pr\'ec\'edente, on obtient le sch\'ema suivant repr\'esentant les ordres partiels entre les diviseurs:

\begin{figure}[H]
$$
\setlength{\unitlength}{1cm}
\begin{picture}(10,4)
\includegraphics{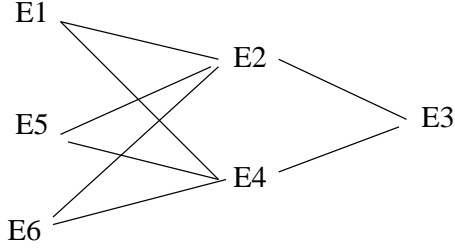}
\end{picture}
$$
\caption{Sch\'ema des ordres partiels}
\end{figure}
(Quand deux $E_i$ sont reli\'es par un trait, avec $Ei$ \`a droite et $E_j$ \`a gauche, cela signifie que $E_i<E_j$, de plus ``l'ordre partiel'' est transitif de la droite vers la gauche ; quand ils ne sont pas reli\'ees, cela veut dire qu il n'y a pas de relation et donc que les familles associ\'ees sont non incluses l'une dans l'autre).
Dans le cas de $E_6$, on a donc les r\'esultats suivants: $\overline{N_1} \not\subset\overline{N_i}$ pour $i \not=1$, $\overline{N_5} \not\subset\overline{N_i}$ pour $i \not=5$, $\overline{N_6} \not\subset\overline{N_i}$ pour $i \not=6$, $\overline{N_2} \not\subset\overline{N_j}$ pour $j=4,3$, $\overline{N_i} \not\subset\overline{N_j}$ pour $j=2,3$, et il reste \`a prouver les autres non inclusions.

\subsection{{\bf Deuxi\`eme application: le cas des $A_n$}}
$A_n$ est la singularit\'e  d'\'equation ($z^{n+1}=xy$). On note les $E_k$ comme ci--dessous:\\

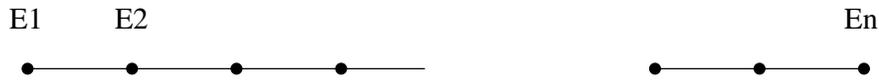
\begin{figure}[H]
$$
\setlength{\unitlength}{1cm}
\begin{picture}(12,1)
\includegraphics{exemple4.pstex}
\end{picture}
$$
\end{figure}

Soit $i<j$.\\
 Soit $f=x$. Alors $ord_{E_i}(f)=i$ et $ord_{E_j}(f)=j$.\\
Donc, par la condition suffisante, $\overline{N_i}\not \subset\overline{N_j}$ .\\
De m\^eme, si $g=y$ , $ord_{E_i}(f)=n-i+1$ et $ord_{E_j}(f)=n-j+1$.\\
Donc $\overline{N_j}\not \subset\overline{N_i}$, toujours en appliquant la condition suffisante.\\
Et ce pour tout $i,j$.\\\\

{\bf La conjecture est prouv\'ee pour les singularit\'es quotients cycliques de type $A_n$}\\\\\\

{\bf Remarque:} Dans le cas des singularit\'es de type $A_n$, il est facile de calculer les familles d'arcs directement et de d\'emontrer qu'elles ne sont pas incluses les unes dans les autres. C'est d 'ailleurs le rpemier exemple trait\'e par Nash dans son preprint [{\bf NASH}]. Nous souhaitons cependant rappeler ce r\'esultat.

Consid\'erons le plongement de $A_n$ dans une vari\'et\'e non singuli\`ere de dimension 3, i.e. l'application surjective $k[[x,y,z]]\rightarrow$${\cal O}_{A_N,0}$. Alors, les arcs sont d\'ecrits par des s\'eries enti\`eres : $(x(t),y(t),z(t))=(a_1t+...,b_1t+...,c_1t+...)$ ; de plus, on ne consid\`ere que les arcs tels que leur transform\'ee stricte intersecte transversalement une des courbes exceptionnelles et n'intersecte pas les autres composantes exceptionnelles; or $z=0$ est l'\'equation locale de chacune des composantes exceptionnelles dans son point g\'en\'erique (on le voit en \'eclatant ``\`a la main'' la singularit\'e d'\'equation ($z^{n+1}=xy$)); cela implique que le premier coefficient de la troisi\`eme s\'erie $c_1$ est non nul sur chacune des vari\'et\'es $\overline{N_i}$. De plus, $ y=\frac{z^{n+1}}{x}$, ce qui implique que les familles de courbes sont:\\
\begin{eqnarray*}
\ N_1={(a_1t+a_2t^{2}+...,b_nt^{n}+...,c_1t+...),  avec\   c_1^{n+1}=b_na_1,...}\\
\ N_2={(a_2t^{2}+...,b_{n-1}t^{n-1}+...,c_1t+...), avec\   c_1^{n+1}=b_{n-1}a_2 ,...}\\
\vdots \\
\ N_n={(a_nt^{n}+...,b_1t+b_2t^{2}+...,c_1t+...), avec\   c_1^{n+1}=b_1a_n,...} 
\end{eqnarray*}\\
Les adh\'erences de ces familles de courbes formelles (pour la topologie de Zariski) sont non incluses les unes dans les autres :\\
Supposons qu'il existe $i,j\in\Nn$ avec $i<j\leq n$ tels que $ \overline{N_i}(A_n)\subset\overline{N_j}(A_n) $.\\
Soit $f\in N_i(A_n)$, $f=(a_it^{i}+...,b_{n-i+1}t^{n-i+1}+...,c_1t+...)$.\\
Alors, $f\in \overline{N_j}(A_n)$, i.e. tout ouvert voisinage de f dans l'espace des arcs intersecte ${N_j}(A_n)$.\\
Pour trouver une contradiction, il nous faut trouver un voisinage $V_f$ de f tel qu'il n'intersecte pas $N_j(A_n)$. Soit $V_f=\{(x(t),y(t),z(t))\in \mathcal{H}/a_i\not=0\}$ ; c'est un voisinage de f ; pour tout $\phi \in V_f$, le i\`eme coefficient de la premi\`ere s\'erie enti\`ere est non nul, alors que tous les \'el\'ements de $N_j(A_n)$ ont la i\`eme coordonn\'ee nulle (ici $i<j$).\\
Donc $V_f$ ne rencontre pas $N_j(A_n)$. $=><=$\\
Donc pour tout $i,j$ on a $\overline{N_i}(A_n)\not\subset\overline{N_j}(A_n)$, i.e. on a $n$ familles d'arcs, dites de Nash .\\\\

\section {{\bf Deuxi\`eme condition suffisante et r\'esolution pour les singularit\'es minimales}} 

\subsection{{\bf Deuxi\`eme condition suffisante}}

Soient $(X,x)$ et $( X_0,x_0)$ deux surfaces alg\'ebriques avec une singularit\'e rationnelles, $p:(X',\{E_i\}_{i\in {\delta} }) \rightarrow (X,x)$ et $p_o:(X_0',\{{E_i}^{0}\}_{i\in {\delta}_1}) \rightarrow (X_0,x_0)$ les d\'esingularisations minimales de $(X,x)$ et $(X_0,x_0)$. Soit $\pi:(X,x)\rightarrow (X_0,x_0)$ un {\bf morphisme dominant, birationnel} entre les deux surfaces.\\
L'application $\pi':(X',\{E_i\}_{i\in {\delta} }) \rightarrow   (X_0',\{{E_i}^{0}\}_{i\in {\delta}_1})$ induite par $\pi$, est birationnelle et donc d\'ecomposable en une suite d'\'eclatements, par cons\'equent on peut associer \`a chaque $\{{E_i}^{0}\}$ une courbe exceptionnelle $E_i$, i.e. $\delta={\delta}_1\cup{\delta}_2$ et ${E_i}^{0}=\pi(E_i)$ pour tout $i\in{\delta}_1$, les autres composantes exceptionnelles se contractant en les points $(a_{i,l})\in {E_{i}}^{0} \ i\in {\delta}_1$ et $l\in {\delta}_2$ .\\
Le morphisme $\pi$ induit sur l'espace des arcs un morphisme continu ${\pi}_{*}$.
Alors on a la proposition suivante:\\
\begin{proposition}
 Si $ \overline{N_i}((X_0,x_0))\not\subset\overline{N_j}((X_0,x_0))$ alors $ \overline{N_i}((X,x))\not\subset\overline{N_j}((X,x)) $ pour $i,j \in{\delta}_1$.
\end{proposition}
\begin{proof}
On a besoin des deux lemmes suivants pour prouver la proposition:\\

\begin{lemme}
 ${\pi}_{*} :N_i((X,x))\rightarrow N_i((X_0,x_0))$ est elle--m\^eme dominante.
\end{lemme}

\begin{proof}

{\bf Notation:} soient $M_i$ l'ensemble des arcs transverses \`a $E_i$, ${M_i}^{0}$ l'ensemble des arcs transverses \`a ${E_i}^{0}$ dans $X'$ et $X_0'$ ; $N_i(X,x)=p(M_i)$ et $N_i(X_0,x_0)=p_0({M_i}^{0})$ ($p$ et $p_0$ sont des isomorphismes hors de l'ensemble exceptionnel). \\

 Nous allons montrer que ${{\pi}_{*}(M_i)}$ est dense dans ${M_i}^{0}$, ce qui impliquera que $\overline {{\pi}_{*}(M_i)}=\overline {{M_i}^{0}}$ et que
 ${\pi}_{*}(N_i((X,x)))$ est dense dans $N_i((X_0,x_0))$.\\
On peut \'ecrire explicitement  ${{\pi}_{*}(M_i)}$:
\begin{eqnarray}
\  {{\pi}_{*}(M_i)}=\{\phi\  arcs / \phi\ transverse\ a\ {E_i}^{0}\ et\ \phi(0)\not = a_{i,l}\ pour\ tout\ l\}.
\end{eqnarray}

Il nous reste maintenant \`a d\'emontrer le lemme suivant:\\

\begin{lemme}
Soient L une surface lisse, $E$ une courbe rationnelle sur L et ${E}^{\#}={E\backslash F}$  o\`u $F=\{a_l\} $ un ensemble fini de points de E.\\
Soient $N=\{\phi\ arcs / \phi\ transverse\ a\ {E}\}$ et $M=\{\phi\ arcs / \phi\ transverse\ a\ {E}^{\#}\}$.\\
Alors $M$ est dense dans $N$.
\end{lemme}

\begin{proof}
On a $M\subset N$ car $M=\{ \phi \in N/ \phi(0) \not \in F \}$ .\\
Donc $M$ est un ouvert de Zariski de $N$, non vide. \\
D'ou $M$ est dense dans $N$.

\end{proof}

Par le lemme pr\'ec\'edent on sait que ${{\pi}_{*}(M_i)}$ est dense dans ${M_i}^{0}$, ce qui implique que ${\pi}_{*}(N_i((X,x)))$ est dense dans $N_i((X_0,x_0))$ .\\

\end{proof}
\bigskip

Par le lemme $3.1$ on sait que ${\pi}_{*}(N_i((X,x)))$ est dense dans $N_i((X_0,x_0))$. Donc il existe $f\in{\pi}_{*}(N_i((X,x)))$ telle qu'il existe un voisinage $V_f$ de $f$ qui  n'intersecte pas $N_j((X_0,x_0))$ (car $ \overline{N_i}((X_0,x_0))\not\subset\overline{N_j}((X_0,x_0)) $ ).
Alors si on pose $g={{\pi}_*}^{-1}(f) \in N_i((X,x))$, et $V_g={{\pi}_*}^{-1}(V_f)$ un voisinage ouvert de g (${\pi}_{*}$ est continue), $V_g$ n'intersecte pas $N_j((X,x))$, ce qui \'equivaut \`a dire que $ \overline{N_i}((X,x))\not\subset\overline{N_j}((X,x)) $.\\\\
\end{proof}

\bigskip

{\bf Exemple d'application}:nous donnerons dans la section suivante un exemple d'application simple (car les singularit\'es minimales et les singularit\'es de type $A_n$ v\'erifient les hypoth\`eses de la proposition), mais nous souhaitons ici donner un autre exemple:\\

\begin{corollaire} 
Supposons que la conjecture est v\'erifi\'ee pour les points doubles rationnels de type $D_n$, $n\geq 4$; alors si $(S,s)$ est un singularit\'e dont le graphe dual de la d\'esingularisation minimale a la m\^eme configuration que celui d'un  $D_n$ (mais pas les m\^emes poids), alors la conjecture est aussi vraie pour $(S,s)$.
\end{corollaire}
\bigskip
{\bf Remarque}:une esquisse de la preuve pour $D_5$ est propos\'ee dans [{\bf REG}]; on peut donc appliquer le corollaire et prouver la conjecture, par exemple, pour la singularit\'e dont le graphe est le suivant:\\

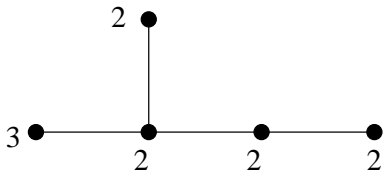
\begin{figure}[H]
$$
\setlength{\unitlength}{1cm}
\begin{picture}(12,4)
\includegraphics{exemple8.pstex}
\end{picture}
$$
\end{figure}

et qui nous promettait de longs calculs pour prouver la conjecture ``\`a la main''.

\subsection{{\bf R\'esolution pour les singularit\'es minimales de surface}} 

\hspace{2cm}Soient $(S,s)$ une singularit\'e minimale, $G$ le graphe dual  de sa d\'esingularisation minimale, $n$ le nombre de sommets de ce graphe. Nous savons d\'ej\`a que le nombre de composantes irr\'eductibles de famille des courbes formelles passant par $(S,s)$ est major\'ee par $n$ (cf. [{\bf REG}] ou [{\bf LEJ}], appendice III), et nous voulons prouver qu'il y a \'egalit\'e.\\

Soient $x$, $y$ deux  sommets de $G$, $N_x$, $N_y$ les familles de courbes associ\'ees \`a $x$ et $y$ respectivement (selon les notations d'Ana Reguera). On veut montrer que les adh\'erences de ces deux familles ne sont pas incluses l'une dans l'autre. La strat\'egie consiste \`a montrer que l'on peut d\'ecomposer en joint birationnel la singularit\'e en singularit\'es de type $A_n$ de telle fa\c{c}on que $x$ et $y$ appartiennent  au m\^eme $A_m$, puis \`a se servir du fait que $ \overline{N_x}(A_m)\not\subset\overline{N_y}(A_m) $ (et inversement) implique que $ \overline{N_x}((S,s))\not\subset\overline{N_y}((S,s)) $ (et inversement) en utlisant la proposition-condition pr\'ec\'edente.\\\\\\ 
 {\bf Notation}:
\begin{itemize}
   \item soient $ z_1 ,z_2,...z_k $ les extr\'emit\'es de $G$ ;
   \item $G_i$ d\'esignera le graphe d'une singularit\'e $A_i$.
\end{itemize}

\bigskip

Par les propri\'et\'es du graphe, on peut joindre de fa\c{c}on unique les deux sommets $x$ et $y$, par un sous-graphe ``en  bambou'' et prolonger celui-ci de telle facon que:\\
\begin{enumerate}
  \item on obtient un nouveau sous--graphe $G$' ``en bambou'' avec $m$ sommets;
  \item les extr\'emit\'es de G' sont par exemple $z_1$ et $z_2$ (i.e. des extr\'emit\'es de $G$ lui- m\^eme, quitte \`a  r\'e\'enum\'erer celles--l\`a) .
\end{enumerate}

Alors le choix de ce sous graphe d\'etermine la d\'ecomposition en joint birationnel de $G$ en graphes $G_i$, en effet :
fixons $z_1$ une extr\'emit\'e ; si $w(z_1)>\gamma(z_1)+1 $, on attache $w(z_1)-\gamma(z_1)-1 $ sommets de poids 1 \`a $z_1$, sinon on ne le touche pas. A chacun des autres sommets z, on attache $r(z)=w(z)-\gamma(z)$ sommets de poids 1
(ici $w(z)$ d\'esigne moins le nombre d'autointersection de la courbe exceptionnelle du graphe de la r\'esolution minimale de $(S,s)$ associ\'ee au sommet $z$ du graphe dual, et $\gamma(z)$ d\'esigne le nombre d'ar\^etes partant de $z$ dans le graphe dual).\\

On a alors plong\'e $G$ dans un graphe non singulier et ce plongement donne la d\'ecomposition de $G$ en $G_i$ : chaque $A_i$ correspond \`a un sous--graphe d'extr\'emit\'es $z_1$ et un sommet de poids $1$ ; il y a autant de $A_i$ que de sommets de poids $1$, et en particulier {\bf un de ces sous graphes correspondant \`a un $A_m$ contient $x$ et $y$.}\\

\begin{figure}[H]
$$
\setlength{\unitlength}{1cm}
\begin{picture}(10,8)
\includegraphics{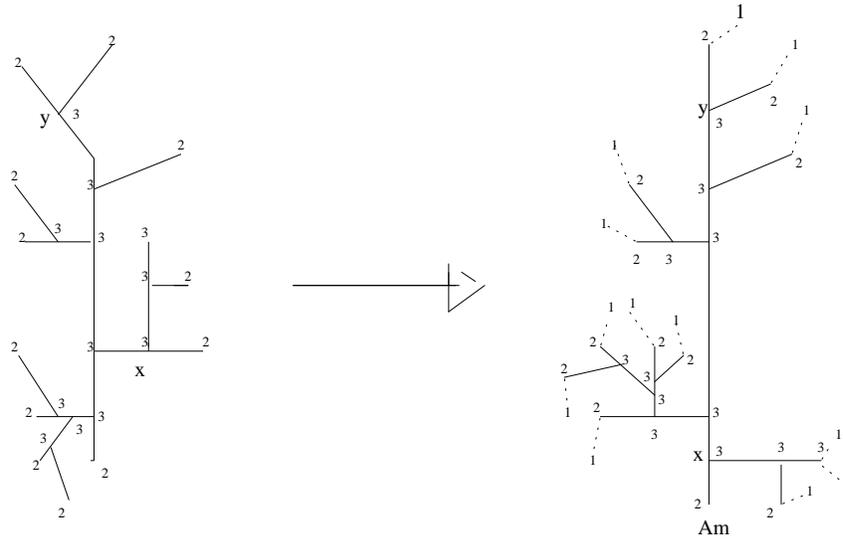}
\end{picture}
$$
\caption{exemple de d\'ecomposition d'un graphe G de singularit\'e minimale en graphes de $A_i$}
\end{figure}

Cette d\'ecomposition de $G$ en $G_i$ correspond au joint birationnel des $A_i$ formant $(S,s)$.On a donc $\pi:(S,s)\rightarrow A_m$ projection qui va de la singularit\'e minimale \`a la singularit\'e $A_m$ ``contenant'' $x$ et $y$; elle est propre, dominante et birationnelle.On peut donc appliquer la proposition pr\'ec\'edente car les singularit\'es de type $A_n$ v\'erifient la conjecture.\\
On a ainsi: $ \overline{N_x}((S,s))\not\subset\overline{N_y}((S,s)) $ et ce, pour tous les $x$, $y$ sommets du graphe dual de r\'esolution de la singularit\'e minimale $(S,s)$

La conjecture est donc r\'esolue pour ces singularit\'es.

\end{document}